\documentstyle{amsppt}
\magnification=1200
\TagsOnRight
\rightheadtext {DUAL APPROACH TO CERTAIN QUESTIONS}
\leftheadtext{B.N.~Khabibullin}
\def\const{\operatorname{const}}
\def\loc{\operatorname{loc}}
\def\Re{\operatorname{Re}}
\def\Spr{\operatorname{Spr}}
\def\supp{\operatorname{supp}}
\define\C{\Bbb C}

\topmatter
\title
DUAL APPROACH TO CERTAIN QUESTIONS
FOR  WEIGHTED SPACES OF HOLOMORPHIC FUNCTIONS.
\endtitle
\author
B.N.~Khabibullin
\endauthor
\abstract
We set new dual problems
for the weighted spaces of
holomorphic functions of one  variable in domains
on $\Bbb C$, namely: nontriviallity of a given space,
description of zero sets,
description of (non-)uniqueness sets,
the representation of meromorphic functions as a quotient
of holomorphic functions from a given space.
\endabstract
\endtopmatter
\document

\head
\S~1. Main results.
\endhead

Let $M$ be a continuous real-valued function in the
domain $G\subset \C$, $0\in G$, and let $H(G)$ be
the space of all holomorphic functions on $G$. We assume
$$
H(G; M)\overset\text{def}\to=
\{ f\in H(G):  \log |f(z)|\leqslant M(z)+\const \, ,\;  z\in G \} \, .
$$

Let $\Lambda =\{ {\lambda}_n\}$ be a sequence
of  complex points in $G$,
${\lambda}_n\to \partial G$, where $\partial G$  is the boundary of
the domain $G$.

We formulate a certain problems for the space $H(G; M)$.
\subhead
1. The problem of the nontriviality of $H(G; M)$
\endsubhead
When does  $f\in H(G; M)$ exist, $f\ne 0$ ?
\medpagebreak

\subhead
2. The problem of the description of zero sets for $H(G; M)$
\endsubhead
We say that $\Lambda$ is a zero set of a function $f\in H(G)$
if the multiplicity of the zero of $f$ at each point $\lambda \in G$
coincides with the number of occurrences of $\lambda$ in the
sequence $\Lambda$. We call a sequence $\Lambda$ a zero set for $H(G; M)$
if there is an $f\in H(G; M)$ with the zero set $\Lambda$.
When is $\Lambda$ a zero set for $H(G; M)$ ?
\subhead
3. The problem of the description of sets of uniqueness
for $H(G; M)$
\endsubhead
We say that $f$ vanishes on $\Lambda$,  and write $f(\Lambda )=0$,
if the multiplicity of the zero of $f$ at each point $\lambda \in G$
is not less than the number of repetitions of $\lambda$ in the
sequence $\Lambda$. The sequence $\Lambda$  is  a set
of nonun\-i\-q\-u\-e\-ness
for a linear space $H\subset H(G)$ if there is a
nonzero function $f\in H$ such that $f(\Lambda )=0$; in the opposite
case, $\Lambda$ is  a set of uniqueness. When is $\Lambda$
a  set of (non-)uniqueness for $H(G; M)$ ?
\subhead
4. The problem of the representation of meromorphic functions
\endsubhead
Let $f$ be a meromorphic in $G$. The function $f$ can
be represented as the ratio $f=g_0/h_0$ of two holomorphic functions
without common zeros. When can $f$ be represented as the ratio
$f=g/h$ of  two functions $g, h \in H(G; M)$ (possibly
with  special conditions on the behavior of $g$ and $h$)?
It is necessary to find the answer in  terms of the function
$$
u_f=\max \{ \log |g_0|, \log |h_0|\} \, .
\tag1.1
$$
Consideration of this function $u_f$ is useful, since, in particular,
$u_f$ is used to define the Nevanlinna characteristic of $f$.
\medpagebreak

The following results give dual formulation for the stated problems.

By $\Cal H(G)$ we denote the space of function harmonic  in $G$.

By $\Cal {SH}(G)$ we denote the class of function subharmonic  in $G$,
$-\infty \notin \Cal {SH}(G)$.

\definition{Definition 1.1} A positive
Radon measure $\mu$ with compact
support in $G$ will be called a representing measure
on $G$ (or with respect to $\Cal H(G)$) if
$$
h(0)= \int h \, d\mu
$$
for any $h\in \Cal H(G)$.

The measure $\mu$ will be called
a Jensen measure on $G$ (or with respect to $\Cal {SH}(G)$) if
$$
u(0)\leq  \int u \, d\mu
$$
for any $u\in \Cal {SH}(G)$.
\enddefinition

Evidently, every Jensen measure is a representing measure.

\definition{Definition 1.2} We call a function $V$ that is subharmonic  in
$G \backslash \{ 0\}$ a representing function on $G$ (or for $\Cal H(G)$)
if is satisfies the following two conditions:

1) There is  a compact $K\subset G$ such that
$V(z)\equiv 0$ for $z\in G\backslash K$;

2) $V(z)\leq -\log |z| +O(1)$ as $z\to 0$. \newline

If $V$  satisfies also the condition

3)  $V(z)\geq 0$ for $z\in G \backslash \{ 0\}$,

then we call $V$ a  Jensen function on $G$ (or for $\Cal{SH}(G)$).
\enddefinition

Evidently, every Jensen function is a representing function.

We use the following notation:

1) For a bounded domain $G$
$$
l_G(z)\overset\text{def}\to= \log \frac1{\rho (z, \partial G)}\, ,
$$

2) For $G=\Bbb C$
$$
l_G(z)\overset\text{def}\to= \log (2+|z|)\, ,
$$

3) For a unbounded domain $G\neq \Bbb C$
$$
l_G(z)\overset\text{def}\to=\max \bigl\{ \log (2+|z|),
\log \frac1{\rho (z, \partial G)}\bigr\} \, .
$$

Let $\sigma (z)$ be a function on $G$ and
$$
0<\sigma (z)<\min \{ \rho (z, \partial G), |z|\}, \quad z\in G,
\tag1.2
$$
where  $\rho (z, \partial G)$
is the distance from $z$ to the boundary $\partial G$ of $G$.

Let $w$ be a function with values in $[-\infty, +\infty ]$.
We put
$$
w^{(\sigma )}(z)\overset\text{def}\to=
\sup\limits_{|\zeta |\le \sigma (z)} w(z+\zeta )\, .
$$
 \proclaim{Theorem 1.1} If $H(G; M)$ is nontrivial, then there exists
 a constant $c\in \Bbb R$ such that
$$
c\leq \int M\, d\mu
\tag1.3
$$
for any Jensen measure $\mu$ on $G$.

Conversely, if  {\rm (1.3)} holds for any Jensen measure $\mu$ on $G$, then
for  the function
$$
\sigma (z)=\min \{ 1, \rho (z, \partial G)/2 \}
\tag1.4
$$
the space $H(G; M^{(\sigma )}+4l_G)$ is nontrivial.
\endproclaim

Let $g_{\Lambda}$ be a function with the zero set $\Lambda$ and
let $M$ be a subharmonic function with mass distribution ${\nu}_{_M}$.

\proclaim{Theorem 1.2} Let $G$ be a simply connected  domain
containing 0. The following three statements   are equivalent:
\roster
\item The sequence $\Lambda$ is a zero set for $H(G; M)$.
\item There is a constant $C$ such that
$$
\sum\limits_n V({\lambda}_n) \leq \int V\, d{\nu}_{_M} +C
\tag1.5
$$
for all  representing  functions $V$ on $G$.
\item There is a constant $C$ such that
$$
\int \log |g_{\Lambda}| \, d\mu  \leq \int M\, d\mu  +C
\tag1.6
$$
for  every representing measure  $\mu$ on $G$.
\endroster
\endproclaim

\proclaim{Theorem 1.3} Let $G$ be a simply connected domain
containing 0. If the sequence $\Lambda$ is  a set of nonuniqueness
for $H(G; M)$, then  there is a constant $C$
such that {\rm (1.5)} holds for any   Jensen function $V$ and
{\rm (1.6)} holds for any  Jensen measure $\mu$ on $G$.

Conversely, if {\rm (1.5)} or {\rm (1.6)} holds for any Jensen function
$V$ on $G$ and   for any  Jensen measure  $\mu$ on $G$,
then for the function $\sigma$ defined by {\rm (1.4)}
the sequence $\Lambda$ is  a set of nonun\-i\-q\-u\-e\-ness
for $H(G; M^{(\sigma )}+4l_G)$.
\endproclaim

In the problems under consideration, the
representing and Jensen functions
play the role of test functions, analogous to the role of the basic
nonnegative functions in the definition of the relation $\leq$
in the space of distributions or in a measure space.

\proclaim{Theorem 1.4} Let $f$ be  meromorphic in $G$.
The function $f$ can be represented in the form $f=g/h$, where $g$ and $h$
are  functions of the class $H(G; M)$ without common zeros,
if and only if there exists a constant $C$ such that
$$
\int u_f\, d\mu \leq \int M\, d\mu +C \,
\tag1.7
$$
for every representing measure  $\mu$ on $G$.

If {\rm (1.7)} holds for any  Jensen measure $\mu$ on a
simply connected domain $G$,
then for the function $\sigma$ defined by  {\rm (1.4)}
the function $f$ can be represented in the form $f=g/h$, where $g$ and $h$
are  functions of the class $H(G; M^{(\sigma )}+4l_G)$.
\endproclaim

These theorems reduce the stated problems to uniform estimates of
the integrals of  certain special classes of functions or
measures.

Theorems 1.1--1.4 are based on a general result,
whose statement follows shortly.

\definition{Definition 1.3}
Let $H$ be a convex cone in $\Cal {SH}(G)$, $0\in G$. A positive
Radon measure $\mu$ with compact
support in $G$ will be called a Jensen measure
with respect to the cone $H$ if
$$
h(0)\leq  \int h \, d\mu
$$
for any $h\in H$.

The class of all Jensen measures with respect to $H$ is
denoted by $\Cal M_H$.
\enddefinition

By   $\, L^1_{\loc}(G)\,$ we denote the space of locally integrable
functions on $G$. A sequence $\{  w_n\} \subset  L^1_{\loc}(G)\,$
is convergent in $L^1_{\loc}(G)\,$  if there exists a function
$w\in L^1_{\loc}(G)\,$ such that
$$
\lim\limits_{n\to \infty} \int\limits_K |w_n-w|\, dm \to 0
$$
for any compact $K\subset G$, where $m$ is the Lebesgue measure on $G$.

\proclaim{Main Theorem} Let $M$ be a continuous real-valued function
on a domain $G\subset \Bbb C$, $0\in G$, $u\in \Cal {SH}(G)$,
and $H$ a convex cone in $\Cal {SH}(G)$ containing the functions
equal to $0$ and $-1$.
If $H$ is sequentially closed in  $\, L^1_{\loc}(G)\,$
and if, for any $h\in H$,  there exists a decreasing sequence
of continuous functions $\{ h_n\} \subset H$ such that
$\; h_n(z)\to h(z)\, $ for $z\in G$ as $\, n\to +\infty$,
then the following statements are equivalent:

1)  There is a $v\in H$ such that
$$
v(0)\neq -\infty \quad \text{ and } \quad
u(z)+v(z)\leq M(z) \, , \quad z\in G \, ;
\tag1.8
$$

2)  There exists a constant $C$ such that
$$
\int u\, d\mu \leqslant \int M \, d\mu + C  \, , \quad \mu \in \Cal M_H \, .
\tag1.9
$$
\endproclaim

The cones $H=\Cal{H} (G)$ and $H=\Cal{SH} (G)$  satisfy of the
conditions of the main theorem.
See for example \cite{1, Theorem 4.1.9}
(sequential closure in $\, L^1_{\loc}(G)\,$ ),
\cite{2, \S~2.1}.

Results near  to the main theorem  were proved in \cite{3, III B, C},
in \cite{4, Theorem on least majorant}
and in \cite{5, Main Theorem}
for different conditions on the cone $H$.

For $G=\C$ the theorems 1.2--1.4 were proved in \cite{5} and \cite{6}.

\head
\S~2. The proof of the main theorem.
\endhead

 Let (1.8) be satisfied.  If we integrate with respect to
an arbitrary Jensen measure $\mu$, we obtain
$$
v(0)\leqslant \int v\, d\mu \leqslant \int (M-u )\, d\mu \, .
$$
This yields the relation  (1.9) for $C=v(0)$.

  Conversly, let (1.9) hold. Then the function
$$
F=M-u
\tag2.1
$$
satisfies
$$
-C\leq \int F\, d\mu \, , \quad \mu \in \Cal M_H\, .
\tag2.2
$$

Let a continuous function $\sigma (z)$ satisfy the condition (1.2).

We denote by $m_z^{(\sigma )}$ the  unit mass uniformly distributed
over the circle $\{ \zeta : |\zeta -z|<\sigma (z) \}$.
If $v$ is a subharmonic function on $G$, then, by the definition of
subharmonicity,
$$
v(z)\leq \int v\, dm_z^{(\sigma )}\, , \quad z\in G, \quad v\in \Cal {SH}(G)
\tag2.3
$$
Consider the function
$$
F_{\sigma  }(z)\overset\text{def}\to=
\int F\, dm_z^{(\sigma )}\, , \quad z\in G.
\tag2.4
$$
The function $F_{\sigma  }$ is continuous since F in (2.1) is
locally Lebesgue integrable.

It follows from the continuity of $M$ that the function $\sigma$
can be chosen   so that
$$
M_{\sigma}(z)\leq M(z)+1 \, , \quad z\in G\, .
\tag2.5
$$

Let $\mu \in \Cal M_H$. Denote by ${\mu}_{\sigma}$ the integral of the
collection of the measures $m_z^{(\sigma )}, z\in G$, with respect to
$\mu$ \  \cite{7}. The Radon measure ${\mu}_{\sigma}$ operate
for $f\in C(G)$ as
$$
\int f \, d{\mu }^{(\sigma )} =
\int \int f\, dm_z^{(\sigma )}\, d\mu (z) \; .
\tag2.6
$$
\proclaim{Lemma 2.1}If $\mu \in \Cal M_H$, then
${\mu }^{(\sigma )} \in \Cal M_H$.
\endproclaim
\demo{Proof} Let $v\in H$. Then it follows from (2.3) and definition 1.3
that
$$
 v(0) \leq \int v\, d\mu \leq
\int \int v\, dm_z^{(\sigma )}\, d\mu (z)
=\int v\, d{\mu}_z^{(\sigma )} \, .
$$
The lemma is proved.
\enddemo

\proclaim{Lemma 2.2} If (2.2) holds, there exists a constant $c>-\infty$
such that the function $F_{\sigma  }$ satisfies
$$
c\,  < \, \int F_{\sigma }\, d\mu \, , \quad  \mu \in \Cal M_H\, .
\tag2.7
$$
\endproclaim
\demo{Proof} By lemma 2.1  it follows from (2.2), (2.4) and (2.6)
that, for any $c<-C$,
$$
c< -C\leq \int F\, d{\mu}^{(\sigma )}
=\int \int F\, dm_z^{(\sigma )}\, d\mu = \int F_{\sigma }\, d{\mu}
$$
for any $\mu \in \Cal M_H$.
The lemma is proved.
\enddemo
\proclaim{Lemma 2.3} Let $H_c=H\cap C(G)$ be the subcone
consisting of all continuous
functions  of the cone $H$.

Then $\Cal M_{H_c}=\Cal M_H$.
\endproclaim
\demo{Proof}  Evidently  $\Cal M_H \subset  \Cal M_{H_c}$. Conversely,
suppose that $h\in H$ and $\mu \in \Cal M_{H_c}$. Then
there exists a decreasing sequence
of continuous functions $\{ h_n\} \subset H$ such that
$\; h_n(z)\to h(z)\, $ as $\, n\to +\infty$ for each  $z\in G$. Consequently,
$$
h(0)\leq \lim_n h_n(0) \leq \lim_n \int h_n \, d\mu
=\int  \lim_n h_n\, d\mu  =\int h\, d\mu .
$$
Therefore $\Cal M_H \supset \Cal M_{H_c}$ and the lemma is proved.
\enddemo

We return to the proof of the theorem.

An exhaustion of the domain $G$ is defined to be a sequence
of subsets $G_n$ that are open and relatively compact in $G$, form a
covering of $G$, and are such that the closure $\overline G_n$
is contained in $G_{n+1}$ for all $n=1, 2,\dots $, i.e.
$G_n\Subset G_{n+1}$.

By $H_n$ we denote the subcone in $C(G_n )$ formed by the restrictions
to $\overline G_n$ of the functions belonging to  $C(G)\cap H$.

The {\it positive germ of the Dirac measure\/} $\delta$
(the unit mass concentrated at zero)  on the cone $H_n$
in the space $C(\overline G_n)$ is
defined \cite{8, III.1.3}, \cite{9, XI, \S~3} 
as the set
$$
{\Spr }_n=\{ \mu \in C^+(\overline G_n): \delta (h)\leq \mu (h),
\; h\in H_n \},
$$
where $C^+(\overline G_n)$ is the cone of positive Radon measures
on $C(\overline G_n)$.
On $C(\overline G_n)$ we define a superlinear functional
$q_n$ \cite{8, p.~299} 
by setting
$$
q_n \, : \, x \to \sup \{ \delta (h) \, : h(z)\leq x(z),
z\in \overline G_n, \; h\in H_n \} , \quad x\in C(\overline G_n).
\tag2.8
$$

The cone $H_n$ is coinitial in $C(\overline G_n)$ \cite{8},
that is, for any $x\in C(\overline G_n)$ there exists an $h\in H_n$
(for example,
a sufficiently small constant) such that $h\leq x$ on $\overline G_n$.
In this case we have the representation \cite{8, III.1.3~(VII)}
$$
q_n(x)=\inf \{ \mu (x)\, :\, \mu \in {\Spr }_n  \},
\quad x\in C(\overline G_n).
\tag2.9
$$
By definition we have $\delta (h)=h(0)$. In other words,
by lemma 2.3 the measures
$\mu \in {\Spr }_n$ are the Jensen measures with respect to $H$,
but under the restriction that the supports of these measures lie in
$\overline G_n$. Thus, (2.9)  can be rewritten for $x=F_{\sigma}$
in the form
$$
q_n(F_{\sigma})=\inf \bigl\{ \int F_{\sigma} \, d\mu \, :
\, \mu \in {\Cal M}_H , \; \supp \mu \in \overline G_n \bigr\}.
$$
By (2.7) the right-hand side is bounded below by a constant $c$,
and so $q_n(F_{\sigma})\geq c$ for all $n$.
Then by definition (2.8) of the functional $q_n$, for some constant
$c$ and for any $n=1, 2,\dots$, there exists a function $h_n\in H_n$
such that
$$
h_n(z)\leq F_{\sigma}(z), \quad z\in \overline G_n \, ;
\qquad h_n(0)\geq c.
\tag2.10
$$

It follows from a  well-known convergence
theorem \cite{1, Theorem 4.1.9 and Proposition 16.1.2} 
that the sequence $\{ h_n\}$ contains a subsequence 
convergent in $\, L_{\loc}^1(G)\,$ to a function $h\leq F_{\sigma}$ on $G$,
$h\not\equiv -\infty$. Since the cone $H$ is
sequentially closed in  $\, L^1_{\loc}(G)\,$
we have $h\in H$. It follows from (2.1) and (2.3)--(2.5)
that
$$
h\leq F_{\sigma}=M_{\sigma}-u_{\sigma} \leq
M+ 1-u  \; \text{ on } G.
$$
We can choose $v=h-1\in H$ and the proof of the theorem is complete.

\head
\S~3. Auxiliary propositions.
\endhead

\proclaim{Proposition 3.1}
Let $M$ and  $u$ be  functions with values in $[-\infty, +\infty ]$
in the simply connected domain $G\subset \C$ and let
$v\in \Cal H (G)$. If $u+v\leq M$
on $G$, then there is a function $h\in H(G)$ without zeros
such that $u+\log |h|\leq M$ on $G$.
\endproclaim

Indeed, if $v\in \Cal H$, then there exists a function $g\in H(G)$, $v=\Re g$.
Set $h=\exp g$.

\proclaim{Proposition 3.2}
Suppose that $M$ is a real-valued continuous function and
 $u$  a subharmonic function on the domain $G\subset \C$.
 If  $v\in \Cal {SH}$ and $u+v\leq M$
on $G$, then for the function $\sigma$ defined by {\rm (1.4)}
there is a function $h\in H(G)$ such that $f\not\equiv 0 $ and
$$
u+\log |h|\leq M^{(\sigma )}+4l_G +\const \; \text{ on }\, G.
$$
\endproclaim
\demo{Proof}  Indeed, by virtue of \cite{6, Lemma 2.2}
for every function $\sigma$ satisfying (1.2) there
exists a function $h\in H(G)$, $h\not\equiv 0 $, such that
$$
u+\log |h| \leq M^{(\sigma )}+4\log (1+|z|+\sigma (z))-2\log \sigma (z),
\quad z\in G.
$$
Choose $\sigma$ according to (1.4). Then for a bounded domain $G$ we have
$$
4\log (1+|z|+\sigma (z))-2\log \sigma (z)\leq \const +
2\log \frac2{\rho (z, \partial G)}\le \const +2l_G(z), \quad z\in G.
$$
For $G=\C$,
$$
4\log (1+|z|+\sigma (z))-2\log \sigma (z)\leq 4\log (1+|z|+1)
\leq 4l_G(z), \quad z\in \C ,
$$
and for an unbounded domain $G$
$$
\multline
4\log (1+|z|+\sigma (z))-2\log \sigma (z) \\
\le 4\log (2+|z|)
+2\log \frac1{\rho (z, \partial G)}+2\log 2 \le
4l_G(z)+\const , \quad z\in G.
\endmultline
$$
The proposition is proved.
\enddemo

Now we establish the interconnection between the Jensen
(re\-s\-p\-e\-c\-t\-i\-vely  re\-p\-re\-s\-e\-n\-ting)  functions
and  the Jensen (respectively  representing)   measures.

Let $\mu$ is a representing measure. We use a potential of the measure $\mu$
equal to the function
(see \cite{6})
$$
V_{\mu}(\zeta )=\int \log |z-\zeta |\, d\mu (z)-\log |\zeta | \, ,
\quad \zeta \in \Bbb C \backslash \{ 0\} \, .
$$

The following proposition can be considered as a Poisson-Jensen
formula.
\proclaim{Proposition 3.3 \cite{6, Proposition 1.4}}
Let $G\subset \C$ be a simply connected  domain containing 0
and let $\mu$ be a representing measure
on G. For every function $u_{\nu}$
subharmonic in $G$ having a mass distribution $\nu$ and such that
$u_{\nu}(0)\ne -\infty$,
$$
\int u_{\nu}\, d\mu =\int V_{\mu}\, d\nu +u_{\nu}(0)\, .
$$
\endproclaim

Denote by $\Cal M_{\Cal H}(G)$, $\Cal M_{\Cal {SH}}(G)$, $\Cal J_{\Cal H}(G)$
and $\Cal J_{\Cal {SH}}(G)$ respectively the classes of
representing measures, of  Jensen measures, of the
representing functions and of the Jensen functions on $G$.

\proclaim{Proposition 3.4 \cite{6, {\rm see proof of
the Proposition 1.5 and 4.1}}}
Let $G\subset \C$ be a simply connected  domain containing 0.
The mapping $P$ defined by the condition
$$
P(\mu )=V_{\mu}, \quad \mu \in \Cal M_{\Cal H}(G)
\quad (\text{ or } \mu \in  \Cal M_{\Cal {SH}}(G)\, )
$$
is a bijection  from the class $\Cal M_{\Cal H}(G)$  onto the class
$\Cal J_{\Cal H}(G)$ (respectively from the class $\Cal M_{\Cal {SH}}(G)$
onto the class $\Cal J_{\Cal {SH}}(G)$).
\endproclaim
For $u\in \Cal {SH}(G)$, we denote the corresponding mass distribution
by ${\nu}_u$.
\proclaim{Proposition 3.5}
Let $G\subset \C$ be a simply connected  domain, $0\in G$.
Let $\Lambda =\{ {\lambda}_n\}$ be a sequence of complex numbers, $0\notin
\Lambda$,  let $g_{\Lambda}$ be a holomorphic function on $G$
with the set $\Lambda$
of zeros, and let $M$ be a function  subharmonic in $G$
with $M(0)\ne -\infty$.The following two statments are equivalent:
$$
\sup\limits_{\mu \in \Cal M_{\Cal H}(G)}
\biggl(\int \log |g_{\Lambda}|\, d\mu -\int M\, d\mu \biggr)<+\infty \, ,
\tag4.1
$$
$$
\sup\limits_{V\in \Cal J_{\Cal H}(G)}
\biggl(\sum\limits_n V({\lambda}_n)-\int V\, d{\nu}_{_M} \biggr)<+\infty \, ,
\tag4.2
$$
and the following two statments are also equivalent:
$$
\sup\limits_{\mu \in \Cal M_{\Cal {SH}}(G)}
\biggl(\int \log |g_{\Lambda}|\, d\mu -\int M\, d\mu \biggr)<+\infty \, ,
\tag4.3
$$
$$
\sup\limits_{V\in \Cal J_{\Cal {SH}}(G)}
\biggl(\sum\limits_n V({\lambda}_n)-\int V\, d{\nu}_{_M} \biggr)<+\infty \, .
\tag4.4
$$
\endproclaim
\demo{Proof} By the generalized of Poisson-Jensen formula from
Proposition 3.3, we obtain
$$
\multline
\int \log |g_{\Lambda}|\, d\mu -\int M\, d\mu\\
=\int V_{\mu}\, d{\nu}_{\log|g_{\Lambda }|} -\int V_{\mu}\, d{\nu}_{_M}
+(\log |g_{\Lambda }(0)|-M(0))\\
=\biggl(\sum\limits_n V_{\mu}({\lambda}_n)-\int V_{\mu}\, d{\nu}_{_M}\biggr)
+(\log |g_{\Lambda }(0)|-M(0)).
\endmultline
$$
Hence, by  Proposition 3.4, relation (4.1) is equivalent to (4.2)
and relation (4.3) is equivalent to (4.4). The present proposition is
proved.
\enddemo

\head
\S~4.  Proof of the theorems 1.1--1.4.
\endhead

\demo{Proof of theorem {\rm 1.1}} If $H(G; M)$ is nontrivial, then there is
a function $f\in H(G)$ such that
$$
\log |f|\le M+C \text{ on } G\, , \quad f\not\equiv 0\, ,
$$
where $C$ is a constant.
The relation (1.3) for any Jensen measure $\mu$ on $G$
follows from the implication $1)\Rightarrow 2)$
of the main theorem for $u\equiv 0$, $v=\log |f| -C$
and $H=\Cal {SH}(G)$.

Conversely, if  {\rm (1.3)} holds for any Jensen measure $\mu$ on $G$,
it follows from the implication $2)\Rightarrow 1)$
of the main theorem for $u\equiv 0$ and $H=\Cal {SH}(G)$ that there exists
a function $v\in \Cal {SH}(G)$  such that $v(0)\ne -\infty$ and
$v(z)\leq M(z), \;  z\in G$. Then by the proposition 3.2
for the function $\sigma$ defined by {\rm (1.4)}
there is a function $f\in H(G)$ such that $f\not\equiv 0 $ and
$\log |f|\leq M^{(\sigma )}+4l_G +\const \; \text{ on }\, G$.
Theorem 1.1 is proved.
\enddemo
\demo{Proof of theorem {\rm 1.2}} If (1) holds, then there exists
a function $f\in H(G, M)$ such that $f=g_{\Lambda}h$,
where $h\in H(G)$ and $h$ does not vanish on $G$. Since for a constant $C$
we have $\log |g_{\Lambda}|+\log h -C \leq M$, it follows from
implication $1)\Rightarrow 2)$
of the main theorem for $u= \log |g_{\Lambda}|$, $v=\log |h| -C$
and $H=\Cal {H}(G)$ that the relation (1.6) is satisfied for  every
representing measure  $\mu$ on $G$, i.~e. $(1)\Rightarrow (3)$.

By virtue of the equivalence $(4.1)\Leftrightarrow (4.2)$ of  proposition
3.5, we obtain the equivalence $(3)\Leftrightarrow (2)$ of  theorem 1.2.

If (3) holds, then, by virtue of the implication $2)\Rightarrow 1)$
of the main theorem for $u= \log |g_{\Lambda}|$ and $H=\Cal {H}(G)$,
there is a harmonic function $v$ on $G$ such that $\log |g_{\Lambda}|+
v\le M$. By  proposition 3.1, there exists a function $h\in H(G)$
such that $\log |g_{\Lambda}|+\log |h|\le M $
and $h$ does not vanish on $G$. We take $f=g_{\Lambda}h\not\equiv 0$. Then
$f\in H(G, M)$ and the sequence $\Lambda$ is a zero set for $f$.
Therefore, the sequence $\Lambda$ is a zero set for $H(G; M)$.

Theorem 1.2 is proved.
\enddemo

\demo{Proof of theorem {\rm 1.3}} If the sequence $\Lambda$ is a set of
nonuniqueness for $H(G; M)$,  there is  a constant $C$
and a function $f\in H(G, M)$ such that $f=g_{\Lambda}h\not\equiv 0$,
$h\in H(G)$.
If $h$  has a zero of multiplicity $k$ at $0$, the analytic function
$h_a(z)=ah(z)/z^k$ satisfies, provided that the positive constant $a$
is sufficiently small, $\log |g_{\Lambda}|+\log |h_a|\le M$ on $G$
and $h_a(0)\ne 0$.  By  the implication $1)\Rightarrow 2)$
of the main theorem for $u= \log |g_{\Lambda}|$, $v=\log |h_a|$
and $H=\Cal {SH}(G)$, the relation (1.6) holds for  every
Jensen measure  $\mu$ on $G$. By virtue of the proposition 3.5,
fulfillment of (1.6) for every  Jensen measure $\mu$ on $G$
is equivalent to  the fulfilment of (1.5) for any Jensen function $V$ on $G$.

If (1.6) is satisfied for every   Jensen measure $\mu$ on $G$,
then, by the implication $2)\Rightarrow 3)$  for $u= \log |g_{\Lambda}|$
and $H=\Cal {SH}(G)$, there exists a function $v\in \Cal {SH}(G)$
such that $v\not\equiv -\infty$ and  $\log |g_{\Lambda}|+v\le M$
on $G$. Hence, by  proposition 3.2, there is a function $h\in H(G)$
such that $h\not\equiv 0$ and
$$
\log |g_{\Lambda}|+\log |h|\leq M^{(\sigma )}+4l_G +\const \; \text{ on }\, G.
$$
Set $f=g_{\Lambda}h\not\equiv 0$. Then $f(\Lambda )=0$ and $f\in
H(G, M^{(\sigma )}+4l_G)$. Therefore, the sequence  $\Lambda$ is
a set of nonuniqueness for $H(G, M^{(\sigma )}+4l_G)$.

Theorem 1.3 is proved.
\enddemo
\demo{Proof of theorem {\rm 1.4}} Suppose that the function $f$
 can be represented in the form $f=g/h$, where $g$ and $h$
are  functions of the class $H(G; M)$ without common zeros.
Then $f=g_0/h_0=g/h$, where  $g_0$ and $h_0$ are the
two holomorphic function without common zeros figuring in (1.1). Hence
$g_0h=gh_0$ and $h/h_0=g/g_0=q$ is a holomorphic function in $G$
without zeros. Therefore,
$$
\multline
u_f+\log |q|=\max \{ \log |g_0|, \log |h_0|\} +\log |q|=
\max \{ \log |g|, \log |h|\}\\
\le M+\const \text{ on } G.
\endmultline
$$
The function $\log |q|$ is harmonic in $G$. By the implication
$1)\Rightarrow 2)$    of the main theorem for $u=u_f$, $v=\log |q|$
and $H=\Cal {H}(G)$,  we obtain (1.7)
for every representing measure $\mu$ on $G$.

Conversely,
if (1.7) holds for every representing measure $\mu$ on $G$,
then, by the implication  $2)\Rightarrow 1)$    of the main theorem
for $u=u_f$ and $H=\Cal {H}(G)$,  there exists a harmonic function
$v$ such that $u_f+v\le M$ on $G$. By virtue of proposition 3.1,
there is a function $q\in H(G)$ without zeros  such that
$$
\max \{ \log |g_0q|, \log |h_0q| \} =u_f+\log |q|\le M \text{ on } G.
$$
Set $g=g_0q\in H(G, M)$ and $h=h_0q\in H(G, M)$. The representation
$f=g/h$ is one sought.

If (1.7) holds for every Jensen measure $\mu$ on $G$,
then, by the implication  $2)\Rightarrow 1)$    of the main theorem
for $u=u_f$ and $H=\Cal {SH}(G)$,  there exists a subharmonic function
$v$ such that $u_f+v\le M$ on $G$. By virtue of proposition 3.2,
there is a function $q\in H(G)$  such that $g\not\equiv 0$ and
$$
\max \{ \log |g_0q|, \log |h_0q| \} =u_f+\log |q|\le
M^{(\sigma )}+4l_G +\const \; \text{ on }\, G.
$$
Set $g=g_0q\in H(G, M^{(\sigma )}+4l_G)$ and
$h=h_0q\in H(G, M^{(\sigma )}+4l_G)$. The representation
$f=g/h$ is the one sought.

Theorem 1.4 is proved.
\enddemo

\head
\S~5. Examples.
\endhead

We illustrate the application of  theorems 1.2 and 1.4 in the
simplest case of $H(\Bbb D , 0)$, where $\Bbb D$ is the unit disk 
in $\C$ and $M(z)\equiv 0$. The characteristics of the new approach
are already quite clear in this case.
\example{Example 1} The classical Nevanlinna theorem shows
that a sequence $\Lambda =\{ {\lambda}_n\} \subset \Bbb D$
is a zero set for $H(\Bbb D , 0)$ if and only if
$$
\sum_n (1-|{\lambda}_n|) \, <\, +\infty \, .
\tag5.1
$$
We will obtain this result as a corollary of  theorem 1.2.

Without loss of generality we can assume that $0 \not\in \Lambda$.

Let $\Lambda$ be a zero set for $H(\Bbb D , 0)$.
The functions
$$
V_r(\zeta )={\log }^+\frac{r}{|\zeta |}\, ,\quad r\, <\, 1,
\tag5.2
$$
where  ${\log }^+(t)\overset\text{def}\to=\max \{\log t, 0\}$,
are  representing functions  on $\Bbb D$

By the implication $1)\Rightarrow 3)$ of  theorem 1.2, we have
$$
\sum_{|{\lambda}_n|\leqslant r}
 \log \frac{r}{|{\lambda}_n|} \leqslant D\, ,\quad r\, <\, 1,
\tag5.3
$$
where the constant $D$ is independent of $r\, <\, 1$. It is easy to show
that the relation (5.3) is equivalent to (5.1).

Conversely, let (5.1) is satiesfied, i.~e., let (5.3) hold.
Let $V$ be a representing function on $\Bbb D$. Then $V$ is
identically equal to $0$ outside of the disk
$\Bbb D \, (r)=\{ \zeta : |\zeta |\, <\, r\}$ for some $r\, <\, 1$.
The function  $V_r$ defined by (5.2) is harmonic in
$\Bbb D\, (r)\, \backslash \{ 0\}$. Therefore, the function $V-V_r$
is subharmonic in $\Bbb D\, (r)\, \backslash \{ 0\}$ and $V-V_r$ is equal
to $0$ for  $|\zeta |=r$.
By virtue of the condition 2) of definition 1.2, the function $V-V_r$
is bounded from above in a deleted neighborhood  of zero. Consequently,
the point $0$ is a removable singularity for the function  $V-V_r$
and the function $V-V_r$ can be continued as a subharmonic function on the
disk   $\, \Bbb D \, (r)$. By the principle of  maximum,
we have $V-V_r\leqslant 0$
on $\Bbb D \, \backslash \{ 0\}$. Hence, by (5.3), we obtain
$$
\sum_n V({\lambda}_n)  \leqslant
\sum_n V_r({\lambda}_n) =
\sum_{|{\lambda}_n|\leqslant r}
 \log \frac{r}{|{\lambda}_n|} \leqslant D\, ,\quad r\, <\, 1.
$$
By  theorem 1.2, this mean that  $\Lambda$  is a zero set
for $H(\Bbb D , 0)$.
\endexample
\example{Example 2} Let $f$ be a mermorphic function of bounded
characteristic in $\Bbb D$. This mean that
$$
T_f(r)\leq \const \; \text{ for all } \; r\, <\, 1\, ,
\tag5.4
$$
where $T_f(r)$ is the Nevanlinna characteristic of $f$.
For simplicity we assume $f(0)=1$.

The Nevanlinna  characteristic of $f$ can be defined as
$$
T_f(r)=\frac1{2\pi}\int\limits_0^{2\pi } u_f (re^{i \theta })\, d\theta \,,
\tag5.5
$$
where $u_f=\max \{ \log |p_0|, \log |q_0|\}$, $f=p_0/q_0$ and
the holomorphic functions $p_0$ and $q_0$ are without
common zeros, and such that $p_0(0)=q_0(0)=1$.

We will be proving the classical Nevanlinna theorem about the representation
of a meromorphic function of bounded   characteristic in $\Bbb D$
as a corollary of  theorem 1.4.

Let $\mu$ be a representing measure on $\Bbb D$ with support
into a disk of radius $r\, < \, 1$, i.~e., $\supp \mu \subset
\Bbb D\, (r)$.
Let ${\mu}^b$  be the balayage of the measure $\mu$ to the
boundary  $\partial \Bbb D\, (r)$ (see \cite{10}). 
By a property of the balayage, for every subharmonic
function $u$ in $\Bbb D$ we have
$$
u(0)\le \int u\, d\mu \le \int u\, d{\mu}^b\, ,\quad u\in \Cal {SH}(\Bbb D)\, .
\tag5.6
$$
Consequently, the measure ${\mu}^b$ is a representing measure on $\Bbb D$.
The measure ${\mu}^b$ is concentrated  on a circle of radius $r$ with
center at zero. Thus we can consider that $\mu^b$ depends only on
$\theta$ for $z=re^{i\theta}$ and we have the equality
$$
\int\limits_0^{2\pi } e^{ik\theta}\, d{\mu}^b(\theta )=0,\quad
k\in \Bbb Z \backslash \{ 0\} \; ; \qquad
\int\limits_0^{2\pi } 1 \, d{\mu}^b(\theta )=1 \, .
$$
This means that $d{\mu}^b(\theta )=\dsize\frac1{2\pi}\, d\theta$.
Hence, by  (5.6), (5.5) and (5.4), we obtain
$$
\int u_f \, d\mu \leq \const \,
$$
for all representing measures  $\mu$  on $0$.
Therefore,  the condition (1.7) of  theorem 1.4 holds for $M\equiv 0$
and for every representing measure $\mu$, and hence the function $f$
can be represented in the form $f=g/h$,where $g$ and $h$ are
 bounded holomorphic functions in $\Bbb D$  without common zeros.

The Nevanlinna theorem is proved.
\endexample

Applications of  theorem 1.3  for the case $G=\C$
can be found in \cite{6}.

\enddocument